\documentclass[10pt,reqno]{article}
\usepackage{graphicx}
\usepackage{amsmath, amssymb}
\usepackage{mathrsfs}
\newtheorem{theorem}{Theorem}[section]

\newtheorem{lemma}{Lemma}[section]

\newtheorem{definition}{Definition}[section]

\numberwithin{equation}{section}

\newcommand{\Ga}{\alpha}
\newcommand{\Gb}{\beta}
\newcommand{\Ge}{\epsilon}
\newcommand{\Gg}{\gamma}
\newcommand{\Gs}{\sigma}
\newcommand{\Ggs}{s}
\newcommand{\Gl}{\lambda}
\newcommand{\GD}{\Delta}

\newcommand{\GG}{\Gamma}
\newcommand{\GL}{\Lambda}

\newcommand{\BM}{{\bf M}}

\newcommand{\qed}{\ \ensuremath{\square}}
\newcommand{\ds}{\displaystyle}
\newcommand{\pf}{\medskip \noindent {\sl Proof}. ~ }
\newcommand{\p}{\partial}
\renewcommand{\a}{\alpha}

\newcommand{\pd}[2]{\frac {\p #1}{\p #2}}

\newcommand{\eqnref}[1]{(\ref {#1})}
\newcommand{\na}{\nabla}
\newcommand{\Om}{\Omega}
\newcommand{\ep}{\epsilon}

\newcommand{\RR}{\mathbb{R}}

\newcommand{\NN}{\mathbb{N}}

\newcommand{\Acal}{\mathcal{A}}
\newcommand{\Ncal}{\mathcal{N}}
\newcommand{\Scal}{\mathcal{S}}

\newcommand{\Kcal}{\mathcal{K}}

\newcommand{\tScal}{\widetilde{\mathcal{S}}}

\newcommand{\be}{\begin{equation}}
\newcommand{\ee}{\end{equation}}

\def\nm{\noalign{\medskip}}

\begin{document}

\title{Reconstruction of Inhomogeneous Conductivities via the Concept of Generalized Polarization
Tensors\footnotemark[1] \footnotetext[1]{This work was supported
by the ERC Advanced Grant Project MULTIMOD--267184 and the
Ministry of Education, Sciences and Technology of Korea through
NRF grants No. 2009-0090250, 2010-0004091 and 2010-0017532.}
\footnotetext[2]{Corresponding author.}
\footnotetext[3]{Department of Mathematics and Applications, Ecole
Normale Sup\'erieure, 45 Rue d'Ulm, 75005 Paris, France
(habib.ammari@ens.fr, deng@dma.ens.fr).}
\footnotetext[4]{Department of Mathematics, Inha University,
Incheon 402-751, Korea (hbkang@inha.ac.kr, hdlee@inha.ac.kr).} }

\author{Habib Ammari\footnotemark[3]
 \and Youjun Deng\footnotemark[2] \footnotemark[3]
\and Hyeonbae Kang\footnotemark[4]
\and Hyundae Lee\footnotemark[4] }

\date{}

\maketitle

\begin{abstract}
This paper extends the concept of generalized polarization tensors
(GPTs), which was previously defined for inclusions with
homogeneous conductivities, to inhomogeneous conductivity
inclusions.  We begin by giving two slightly different but
equivalent definitions of the GPTs for inhomogeneous inclusions.
We then show that, as in the homogeneous case, the GPTs are the
basic building blocks for the far-field expansion of the voltage
in the presence of the conductivity inclusion. Relating the GPTs
to the Neumann-to-Dirichlet (NtD) map, it follows that the full
knowledge of the GPTs allows unique determination of the
conductivity distribution. Furthermore, we show important
properties of the the GPTs, such as symmetry and positivity,  and
derive bounds satisfied by their harmonic sums. We also compute
the sensitivity of the GPTs with respect to changes in the
conductivity distribution and propose an algorithm for
reconstructing conductivity distributions from their GPTs. This
provides a new strategy for solving the highly nonlinear and
ill-posed inverse conductivity problem. We demonstrate the
viability of the proposed algorithm by preforming a sensitivity
analysis and giving some numerical examples.
\end{abstract}

\noindent {\footnotesize {\bf Mathematics subject classification
(MSC2000):} 35R30, 35C20}


%%%%%%%%%%%%%%%%%%%%%%%%%%%%%%%%%%%%%%%%%%%%%%%%%%%%
\section{Introduction}
%%%%%%%%%%%%%%%%%%%%%%%%%%%%%%%%%%%%%%%%%%%%%%%%%%%%

There are several geometric and physical quantities associated
with shapes such as eigenvalues and capacities \cite{PS51}. The
concept of the generalized polarization tensors (GPTs) is one of
them. The notion appears naturally when we describe the
perturbation of the electrical potential due to the presence of
inclusions whose material parameter (conductivity) is different
from that of the background.

To mathematically introduce the concept of GPTs, we consider the
conductivity problem in $\RR^d$, $d=2,3$:
\begin{equation}\label{condeqn}
\left\{
\begin{array}{ll}
\nabla \cdot (\chi(\RR^d \setminus \overline{\Om}) + k \chi(\Om)) \nabla u =0 & \mbox{in } \RR^d, \\
\nm u(x)-h(x)=O(|x|^{1-d}) & \mbox{as } |x| \rightarrow \infty.
\end{array}
\right.
\end{equation}
Here, $\Om$ is the inclusion embedded in $\RR^d$ with a Lipschitz
boundary, $\chi(\Om)$ (resp. $\chi(\RR^d \setminus
\overline{\Om})$) is the characteristic function of $\Om$ (resp.
$\RR^d \setminus \overline{\Om}$), the positive constant $k$ is
the conductivity of the inclusion which is supposed to be
different from the background conductivity $1$, $h$ is a harmonic
function in $\RR^d$ representing the background electrical
potential, and the solution $u$ to the problem represents the
perturbed electrical potential. The perturbation $u-h$ due to the
presence of the conductivity inclusion $\Omega$ admits the
following asymptotic expansion as $|x| \to \infty$:
 \be\label{uhasymphom}
 u(x)-h(x) = \sum_{|\Ga|, |\Gb| \ge 1} \frac{(-1)^{|\Gb|}}{\Ga ! \Gb !}
 \p^\Ga h(0) M_{\Ga\Gb}(k, \Om) \p^\Gb \Gamma(x),
 \ee
where $\GG$ is the fundamental solution of the Laplacian (see, for
example, \cite{AK03, book2}). The building blocks  $M_{\Ga\Gb}(k,
\Om)$ for the  asymptotic expansion  \eqnref{uhasymphom} are
called the GPTs. Note that the GPTs $M_{\Ga\Gb}(k, \Om)$ can be
reconstructed from the far-field measurements of $u$ by a
least-squares method. A stability analysis of the reconstruction
is provided in \cite{han}. On the other hand, it is shown in
\cite{han2} that in the full-view case, the reconstruction problem
of GPTs from boundary data has the remarkable property that low
order GPTs are not affected by the error caused by  the
instability of higher-orders in the presence of measurement noise.

The GPTs carry geometric information about the inclusion. For
example, the inverse GPT problem holds to be true, namely, the
whole set of GPTs, $\{ M_{\Ga\Gb}(k, \Om):  |\Ga|, |\Gb| \ge 1
\}$, determines $k$ and $\Om$ uniquely \cite{AK02}. The leading
order GPT (called the polarization tensor (PT)), $\{ M_{\Ga\Gb}(k,
\Om):  |\Ga|, |\Gb| = 1 \}$, provides the equivalent ellipse
(ellipsoid) which represents overall property of the inclusion
\cite{AKKL03, BHV03}. Moreover, there are important analytical and
numerical studies which show that finer details of the shape can
be recovered using higher-order GPTs \cite{AKLZ12, AGKLY11}. The
GPTs even carry topology information of the inclusion
\cite{AGKLY11}. It is also worth mentioning that an efficient
algorithm for computing the GPTs is presented in \cite{yves}.

The notion of GPTs appears in various contexts such as asymptotic
models of dilute composites ({\it cf.} \cite{Ko92, Mi01, AKT05}),
low-frequency asymptotics of waves \cite{DK89},  potential theory
related to certain questions arising in hydrodynamics \cite{PS51},
biomedical imaging of small inclusions (see \cite{handbook} and
the references therein), reconstructing small inclusions
\cite{FV89, AKKL03, BHV03}, and shape description \cite{AGKLY11}.
Recently the concept of GPTs finds another promising application
to cloaking and electromagnetic and acoustic invisibility. It is
shown that the near-cloaking effect of  \cite{KSVW08} can be dramatically improved
by using multi-layered structures whose GPTs vanish up to a
certain order \cite{AKLL11}.

As far as we know, the GPTs have been introduced only for
inclusions with homogeneous conductivities or layers with constant
conductivities. It is the purpose of this paper to extend the
notion of GPTs to inclusions with inhomogeneous conductivities and
use this new concept for solving the inverse conductivity problem.
We first introduce the GPTs for inhomogeneous inclusions and show
that exactly the same kind of far-field asymptotic formula as
\eqnref{uhasymphom} holds. We also prove important properties of
the GPTs such as unique determination of Neumann-to-Dirichlet map,
symmetry, and positivity. We then provide a sensitivity analysis
of the GPTs with respect to changes in the conductivity
distribution. We finally propose a minimization algorithm for
reconstructing an inhomogeneous conductivity distribution from its
high-order GPTs. We carry out a resolution and stability analysis
for this reconstruction problem in the linearized case and present
numerical examples to show its viability.

The paper is organized as follows. In section 2 we introduce the
GPTs for inhomogeneous conductivity inclusions and prove that they
are the building blocks of the far-field expansion of the
potential. Section 3 is devoted to the derivation of integral
representations of the GPTs. We also establish a relation between
the GPTs and the NtD map. In section 4 we prove important
properties of symmetry and positivity of the GPTs and obtain
bounds satisfied by their harmonic sums. In section 5 we perform a
sensitivity analysis of the GPTs with respect to the conductivity
distribution. We also show that in the linearized case, high-order
GPTs capture high-frequency oscillations of the conductivity. In
section 6, we present an algorithm for reconstructing
inhomogeneous conductivity distributions from their high-order
GPTs. The algorithm is based on minimizing the discrepancy between
the computed and measured GPTs.

%%%%%%%%%%%%%%%%%%%%%%%%%%%%%%%%%%%%%%%%%%%%%%%%%%%
\section{Contracted GPTs and asymptotic expansions}
%%%%%%%%%%%%%%%%%%%%%%%%%%%%%%%%%%%%%%%%%%%%%%%%%%%

Let $\Gs$ be a bounded measurable function in $\RR^d$, $d=2,3$, such that $\Gs-1$ is compactly supported and
 \be\label{glonetwo}
 \Gl_1 \le \Gs \le \Gl_2
 \ee
for positive constants $\Gl_1$ and $\Gl_2$. For a given harmonic function $h$ in $\RR^d$, we consider the following conductivity problem:
\begin{equation}\label{eq:101}
\left\{
\begin{array}{ll}
\nabla \cdot \sigma \nabla u =0 & \mbox{in } \RR^d, \\
\nm u(x)-h(x)=O(|x|^{1-d}) & \mbox{as } |x| \rightarrow \infty.
\end{array}
\right.
\end{equation}
In this section we derive a full far-field expansion of $(u-h)(x)$
as $|x| \to \infty$. In the course of doing so, the notion of
(contracted) generalized polarization tensors (GPT) appears
naturally.

Let $B$ be a bounded domain in $\RR^d$ with a
$\mathcal{C}^{1,\eta}$-boundary $\partial B$ for some $0<\eta<1$.
We assume that $B$ is such that
 \be\label{suppB}
 \mbox{supp}\, (\Gs-1) \subset B.
 \ee
Suppose that $B$ contains the origin. Let $H^{s}(\p B)$, for $s\in
\RR$,  be the usual $L^2$-Sobolev space and let $H^{s}_0 (\p
B):=\{\phi \in H^{s}(\p B) | \int_{\p B} \phi =0\}$. For $s=0$, we
use the notation $L^2_0(\partial B)$.

The Neumann-to-Dirichlet (NtD) map $\GL_{\Gs} : H^{-1/2}_0(\p B)
\to H^{1/2}_0(\p B)$ is defined to be
 \be
 \GL_{\Gs}[g]:= u|_{\p B},
 \ee
where $u$ is the solution to
\be\label{sigmaeqn}
\left\{
\begin{array}{ll}
\nabla \cdot \Gs\nabla u=0 & \mbox{in }  B, \\
\ds \Gs\frac{\partial u}{\partial \nu} = g & \mbox{on } \ds \partial B \quad (\int_{\p B} u=0)
\end{array}
\right. \ee for $g \in H^{-1/2}_0(\p B)$. The operator $\GL_1$ is
the NtD map when $\Gs\equiv 1$.

Note that \eqnref{eq:101} is equivalent to
\begin{equation}\label{eq:101-1}
\left\{
\begin{array}{ll}
\nabla \cdot\Gs \nabla u=0 & \mbox{in }  B, \\
\GD u =0 & \mbox{in } \RR^d \setminus \overline{B}, \\
\ds \frac{\p u}{\p \nu} \Big|_+=\Gs \frac{\p u}{\p \nu} \Big|_- & \mbox{on }  \p B,\\
u|_+=u|_- & \mbox{on } \partial B,\\
u(x)-h(x)=O(|x|^{1-d}) & \mbox{as } |x| \rightarrow \infty.
\end{array}
\right.
\end{equation}
Here and throughout this paper, the subscripts $\pm$ indicate the
limits from outside and inside $B$, respectively.

Let $\GG(x)$ be the fundamental solution to the Laplacian: \be
\label{gammacond} \Gamma (x) = \begin{cases} \ds \frac{1}{2\pi}
\ln |x|\;, \qquad & d=2, \\ \nm \ds -\frac{1}{4\pi} |x|^{-1}\;,
\qquad & d = 3.
\end{cases}
\ee
If $u$ is the solution to \eqnref{eq:101}, then by Green's formula
we have for $x \in \RR^d \setminus \overline{B}$
\begin{align*}
(u-h)(x) & =  \int_{\p B} \GG(x-y) \frac{\p (u-h)}{\p \nu} \Big
|_{+} (y) d\Ggs_y -
\int_{\p B}  \frac{\p \GG(x-y)}{\p \nu_y} (u-h)|_{+} (y) d\Ggs_y \\
& =  \int_{\p B} \GG(x-y) \frac{\p u}{\p \nu} \Big|_{+} (y)
d\Ggs_y - \int_{\p B}  \frac{\p \GG(x-y)}{\p \nu_y} u|_+ (y)
d\Ggs_y ,
\end{align*}
where the second equality holds since $h$ is harmonic. Let
$g=\Gs\frac{\partial u}{\partial \nu}|_-$. Then we have $u|_{\p
B}=\GL_{\Gs}[g]$ on $\p B$. Thus we get from the transmission
conditions in \eqnref{eq:101-1} that
\begin{align} \label{eqa28}
(u-h)(x) & =  \int_{\p B} \GG(x-y) g(y) d\Ggs_y - \int_{\p B}
\frac{\p \GG(x-y)}{\p \nu_y} \GL_{\Gs}[g](y) d\Ggs_y .
\end{align}
For $x \in \RR^d \setminus \overline{B}$, we have
 $$
 \GL_1 \left( \frac{\p \GG(x - \cdot)}{\p \nu_y} \right) = \GG(x - \cdot) -\frac{1}{|\p B|} \int_{\p B} \GG(x-y) d\Ggs_y \quad \mbox{on } \p B,
 $$
and hence
 \begin{equation} \label{eqb28}
 \int_{\p B}  \frac{\p \GG(x-y)}{\p \nu_y} \GL_{\Gs}[g](y) d\Ggs_y = \int_{\p B} \GG(x-y) \GL_1^{-1} \GL_{\Gs}[g](y) d\Ggs_y .
 \end{equation}
Thus we get from \eqnref{eqa28} and  \eqnref{eqb28} that
 \be\label{uminusH}
 (u-h)(x) = \int_{\p B} \GG(x-y) \GL_1^{-1} (\GL_1- \GL_{\Gs}) [g] (y) d\Ggs_y, \quad x \in \RR^d \setminus \overline{B}.
 \ee
Here we have used the fact that $\GL_1: H^{-1/2}_0(\partial B)
\rightarrow H^{1/2}_0(\partial B)$ is invertible and self-adjoint:
$$
\langle \GL_1[g], f\rangle_{H^{1/2}, H^{-1/2}} = \langle g,
\GL_1[f]\rangle_{H^{1/2}, H^{-1/2}}, \quad \forall\; f,g \in
H^{-1/2}_0(\partial B),
$$
with $\langle , \rangle_{H^{1/2}, H^{-1/2}}$ being the duality
pair between $H^{-1/2}(\partial B)$ and  $H^{1/2}(\partial B)$.

Suppose that $d=2$. For each positive integer $n$, let $u_n^c$ and
$u_n^s$ be the solutions to \eqnref{eq:101} when $h(x)=r^n \cos
n\theta$ and $h(x)=r^n \sin n\theta$, respectively. Let
 \be
 g_n^c := \Gs \frac{\p u_n^c}{\p\nu} \Big|_{-} \quad\mbox{and}\quad g_n^s := \Gs \frac{\p u_n^s}{\p\nu} \Big|_{-} \quad \mbox{on } \p B.
 \ee
Since \eqnref{eq:101} is linear, it follows that if the harmonic
function $h$ admits the expansion
 \be\label{Hexp}
 h(x) = h(0)+\sum_{n=1}^\infty r^n \bigr(a_n^c\cos n\theta + a_n^s\sin n\theta \bigr)
 \ee
with $x=(r\cos\theta, r\sin\theta)$, then we have
 $$
 g:= \Gs\frac{\partial u}{\partial \nu} \big|_{-}=  \sum_{n=1}^\infty \bigr(a_n^c g_n^c + a_n^s g_n^s \bigr),
 $$
and hence
 \be\label{uminusH-2}
 (u-h)(x) = \sum_{n=1}^\infty \int_{\p B} \GG(x-y) \bigr( a_n^c \GL_1^{-1} (\GL_1- \GL_{\Gs}) [g_n^c] (y)
 + a_n^s \GL_1^{-1} (\GL_1- \GL_{\Gs}) [g_n^s] (y) \bigr) d\Ggs_y.
 \ee
Note that $\GG(x-y)$ admits the expansion
 \be\label{Gammaexp}
 \GG(x-y)
 =\sum_{n=1}^{\infty}\frac{-1}{2\pi n}\left[\frac{\cos n\theta_x}{r_x^n}r_y^n\cos n\theta_y
 +\frac{\sin n\theta_x}{r_x^n}r_y^n\sin n\theta_y\right]+ C,
 \ee
where $C$ is a constant, $x=r_x(\cos\theta_x,\sin\theta_x)$ and
$y=r_y(\cos\theta_y,\sin\theta_y)$. Expansion \eqnref{Gammaexp} is
valid if $|x| \to \infty$ and $y \in \p B$. The contracted
generalized polarization tensors are defined as follows (see
\cite{AKLL11}):
 \begin{align}
 \ds M_{mn}^{cc}= M_{mn}^{cc}[\Gs]:= \int_{\p B} r_y^m\cos m\theta_y \, \GL_1^{-1} (\GL_1- \GL_{\Gs})
 [g_n^c] (y) \, d\Ggs_y, \label{defmcc}\\
 \ds M_{mn}^{cs}= M_{mn}^{cs}[\Gs]:=\int_{\p B} r_y^m\cos m\theta_y \, \GL_1^{-1} (\GL_1- \GL_{\Gs}) [g_n^s] (y) \, d\Ggs_y, \\
 \ds M_{mn}^{sc}= M_{mn}^{sc}[\Gs]:= \int_{\p B} r_y^m\sin m\theta_y \, \GL_1^{-1} (\GL_1- \GL_{\Gs}) [g_n^c] (y) \, d\Ggs_y, \\
 \ds M_{mn}^{ss}= M_{mn}^{ss}[\Gs]:= \int_{\p B} r_y^m\sin m\theta_y \, \GL_1^{-1} (\GL_1- \GL_{\Gs}) [g_n^s] (y) \, d\Ggs_y.
 \end{align}
From \eqnref{uminusH-2} and \eqnref{Gammaexp}, we get the
following theorem.
\begin{theorem}
\label{th:expan2} Let $u$ be the solution to \eqnref{eq:101} with
$d=2$. If $h$ admits the expansion \eqnref{Hexp}, then we have
\begin{align}
(u-h)(x) & = -\sum_{m=1}^\infty\frac{\cos m\theta}{2\pi mr^m}\sum_{n=1}^\infty
\bigr(M_{mn}^{cc}a_n^c + M_{mn}^{cs}a_n^s \bigr) \nonumber \\
& \qquad -\sum_{m=1}^\infty\frac{\sin m\theta}{2\pi m
r^m}\sum_{n=1}^\infty \bigr(M_{mn}^{sc}a_n^c + M_{mn}^{ss}a_n^s
\bigr), \label{expan2}
\end{align}
which holds uniformly as $|x| \to \infty$.
\end{theorem}

In three dimensions, we can decompose harmonic functions as
follows:
\begin{align}
h(x)=h(0) + \sum_{n=1}^\infty \sum_{m=-n}^n a_{mn} r^n Y_n^m(\theta, \varphi), \label{hdecom}
\end{align}
where $(r,\theta, \varphi)$ is the spherical coordinate of $x$ and $Y_n^m$ is the spherical harmonic function of degree $n$ and of order $m$. Let
 \be\label{gmn}
 g_{mn}=\Gs \frac{\p u_{mn}}{\p\nu} \Big|_{-} \quad\mbox{on } \partial
 B,
 \ee
where $u_{mn}$ is the solution to \eqnref{eq:101} when $h(x)=r^n
Y_n^m(\theta, \varphi)$. It is well-known (see, for example,
\cite{nedelec}) that
\begin{align}
\Gamma(x-y)= - \sum_{\ell=0}^\infty \sum_{k=-\ell}^{\ell} \frac{1}{2\ell+1}Y_\ell^k(\theta, \varphi) \overline{Y_\ell^k(\theta', \varphi')} \,\frac{r'^n}{r^{n+1}},\label{Gammaexp2}
\end{align}
where $(r,\theta, \varphi)$ and $(r',\theta', \varphi')$ are the
spherical coordinates of $x$ and $y$, respectively. Analogously to
Theorem \ref{th:expan2}, the following result holds.
\begin{theorem}
\label{th:expan3} Let $u$ be the solution to \eqnref{eq:101} with
$d=3$. If $h$ admits the expansion \eqnref{hdecom}, then we have
\be\label{expan3} (u-h)(x) =- \sum_{\ell=1}^\infty
\sum_{k=-\ell}^{\ell}\sum_{n=1}^\infty \sum_{m=-n}^n
\frac{a_{mn}M_{mnk\ell}}{(2\ell+1)r^{n+1}}Y_\ell^k(\theta,
\varphi)\quad\mbox{\rm as } |x| \to \infty, \ee where the GPT
$M_{mnk\ell}=M_{mnk\ell}[\Gs]$ is defined by \be
M_{mnk\ell}:=\int_{\partial B} Y_\ell^k (\theta', \varphi') r'^n
\GL_1^{-1} (\GL_1- \GL_{\Gs})[g_{mn}] (r',\theta', \varphi') \,
d\sigma. \ee
\end{theorem}

We emphasize that the definitions of contracted GPTs do not depend on the choice of $B$ as long as \eqnref{suppB} is satisfied.
This can be seen easily from \eqnref{expan2} and \eqnref{expan3} (see also section \ref{section4}).

%%%%%%%%%%%%%%%%%%%%%%%%%%%%%%%%%%%%%%%%%%%%%%%%%%%
\section{Integral representation of GPTs }
%%%%%%%%%%%%%%%%%%%%%%%%%%%%%%%%%%%%%%%%%%%%%%%%%%%

In this section, we provide another definition of GPTs which is based on integral equation formulations as in
\cite{KS2000, book2}. Proper linear combinations of GPTs defined in this section coincide with the
contracted GPTs defined in the previous section.

Let $N_{\Gs}(x,y)$ be the Neumann function of problem
\eqnref{sigmaeqn}, that is, for each fixed $z \in B$,  $N_{\Gs}(x,
y)$ is the solution to
\begin{equation} \label{revn}
 \ \left \{
\begin{array}{l}
\ds \nabla \cdot \Gs \nabla N(\cdot, z) = -\delta_z(\cdot) \quad
\mbox{in } B , \\ \nm \ds \Gs \nabla N(\cdot, z) \cdot \nu
\big|_{\p B} = \frac{1}{|\partial B |} , \quad \ds \int_{\p B}
N(x,z) \, d\sigma(x) =0 .
\end{array}
\right .\end{equation} Then  the function $u$ defined by
 \be
 u(x)= \Ncal_{B, \Gs} [g](x):=\int_{\p B} N_{\Gs}(x,y) g(y) d\Ggs_y, \quad x \in B
 \ee
is the solution to \eqnref{sigmaeqn}, and hence
 \be\label{GLGsNcal}
 \GL_{\Gs}[g](x) = \Ncal_{B, \Gs} [g](x), \quad x \in \p B.
 \ee

Let $\Scal_B$ be the single layer potential on $\p B$, namely,
 \be\label{eq:sglyer}
 \Scal_B[\phi](x)=\int_{\p B} \GG(x-y) \phi(y)
 d\Ggs_y, \quad x \in \RR^d.
 \ee
Let the boundary integral operator
$\Kcal_B$ (sometimes called the Poincar\'e-Neumann operator) be
defined by
$$
\Kcal_B[\phi](x)= \int_{\p B} \frac{\partial \GG}{\partial
\nu_y}(x-y) \phi(y) d\Ggs_y.
$$
It is well-known that the single layer potential $\Scal_B$
satisfies the trace formula \be \label{eq:trace}
\frac{\p}{\p\nu}\Scal_B[\phi] \Big|_{\pm} = (\pm \frac{1}{2}I+
\Kcal_{B}^*)[\phi] \quad \mbox{on } \p B, \ee where $\Kcal_{B}^*$
is the $L^2$-adjoint of $\Kcal_B$. We recall that $\Gl
I-\Kcal_{B}^*$ is invertible on $L_0^2(\p B)$ if $|\Gl|\geq 1/2$
(see, for example, \cite{Folland76, Ver84, book2}).

Identity \eqnref{uminusH} suggests that the solution $u$ to
\eqnref{eq:101} may be represented as
 \be\label{solrep}
 u(x) = \left \{
 \begin{array}{ll}
 h(x) + \Scal_B [\phi](x), & \ \ x \in \RR^d \setminus B, \\
 \nm \Ncal_{B, \Gs} [\psi](x) + C, &  \ \ x \in B
 \end{array}
 \right.
 \ee
for some densities $\phi$ and $\psi$ on $\p B$, where the
constant $C$ is given by
 \be
 C= \frac{1}{|\p B|} \int_{\p B} \left( h + \Scal_B [\phi] \right) \; d\Ggs.
 \ee
In view of the transmission conditions along $\p B$ in
(\ref{eq:101}), \eqnref{GLGsNcal} and \eqnref{eq:trace}, the pair
of densities $(\phi, \psi)$ should satisfy
 \be\label{inteqn}
 \left \{
 \begin{array}{l}
  \ds -\mathcal{S}_B [\phi] + \frac{1}{|\p B|} \int_{\p B} \Scal_B [\phi] \; d\Ggs + \GL_{\Gs}[\psi] = h - \frac{1}{|\p B|} \int_{\p B} h \, d\Ggs\\
 \nm
 \ds -(\frac{1}{2} I + \mathcal{K}_B^*)[\phi] + \psi = \frac{\p h}{\p \nu}
 \end{array}
 \right.
 \ \ \mbox{on} \ \ \p B.
 \ee

We now prove that the integral equation \eqnref{inteqn} is uniquely solvable. For that, let
 \be
 \tScal_B [\phi] := \mathcal{S}_B [\phi] - \frac{1}{|\p B|} \int_{\p B} \Scal_B [\phi] \; d\Ggs.
 \ee

\begin{lemma}\label{invert}
The operator $\Acal:H^{-1/2}(\p B) \times H_0^{-1/2}(\p B) \to H_0^{1/2}(\p B) \times H^{-1/2}(\p B)$ defined by
\be
 \Acal:= \begin{bmatrix}
 - \tScal_B & \GL_{\Gs} \\
 \nm
 \ds -(\frac{1}{2} I + \mathcal{K}_B^*) & I
 \end{bmatrix}
 \ee
is invertible.
\end{lemma}

As an immediate consequence of Lemma \ref{invert} we obtain the
following theorem.
\begin{theorem}
The solution $u$ to \eqnref{eq:101} can be represented in the form
\eqnref{solrep} where the pair $(\phi, \psi) \in H^{-1/2}(\p B)
\times H_0^{-1/2}(\p B)$ is the solution to
 \be\label{Acalphi2}
 \Acal \begin{bmatrix} \phi \\ \psi \end{bmatrix} = \begin{bmatrix} h- \frac{1}{|\p B|} \int_{\p B} h \; d\Ggs \\
 \nm
 \pd{h}{\nu}|_{\p B} \end{bmatrix}.
 \ee
\end{theorem}

\noindent{\it Proof of Lemma \ref{invert}}. We first recall the
invertibility of
 $\Scal_B: H^{-1/2}(\p B) \to H^{1/2}(\p B)$ in three dimensions
 (see, for instance, \cite{Ver84}). In two dimensions this result is not anymore true. However, using Theorem 2.26 of \cite{book2},
 one can show that in two dimensions there exists a unique $\phi_0 \in L^2(\p B)$ such that
 \be
 \int_{\p B} \phi_0 =1\quad\mbox{and}\quad \tilde{\Scal}_B [\phi_0]=0 \quad \mbox{on } \p B.
 \ee
Then we have
 \begin{equation}\label{phipsieq}
 \Acal\begin{bmatrix} \phi_0 \\ 0 \end{bmatrix} =  \begin{bmatrix} 0 \\ -(\frac{1}{2} I + \mathcal{K}_B^*)[\phi_0] \end{bmatrix},
 \end{equation}
and
 \be\label{kstarone}
 \int_{\p B} \left(\frac{1}{2} I + \mathcal{K}_B^*\right)[\phi_0] d\sigma = \int_{\p B} \phi_0 \left(\frac{1}{2} I + \mathcal{K}_B \right)[1] d\sigma = \int_{\p B} \phi_0 d\sigma = 1 .
 \ee
Therefore, by replacing $\phi$ with $\phi - \phi_0 \int_{\partial
B} \phi$,  it is enough in both the two- and three-dimensional
cases to determine uniquely $(\phi, \psi) \in H_0^{-1/2}(\p B)
\times H_0^{-1/2}(\p B)$ satisfying
 \be\label{Acalphi}
 \Acal \begin{bmatrix} \phi \\ \psi \end{bmatrix} = \begin{bmatrix} f \\ g \end{bmatrix}
 \ee
for $(f,g) \in H_0^{1/2}(\p B) \times H_0^{-1/2}(\p B)$. In fact,
if $(f,g) \in H_0^{1/2}(\p B) \times H^{-1/2}(\p B)$, then let
$C=\frac{1}{|\p B|} \int_{\p B} g$ and let $(\phi, \psi)$ be the
solution to
 $$
 \Acal \begin{bmatrix} \phi \\ \psi \end{bmatrix} = \begin{bmatrix} f \\ g-C (\frac{1}{2} I + \mathcal{K}_B^*)[\phi_0] \end{bmatrix} .
 $$
It then follows from \eqnref{phipsieq} and \eqnref{kstarone} that
 $$
 \Acal \begin{bmatrix} \phi - C \phi_0 \\ \psi \end{bmatrix} = \begin{bmatrix} f \\ g \end{bmatrix} .
 $$

We now show that \eqnref{Acalphi} is uniquely solvable for a given
$(f,g) \in H_0^{1/2}(\p B) \times H_0^{-1/2}(\p B)$. We first
introduce the functional spaces $$H^1_{\mathrm{loc}}(\RR^d):=\{ h
u \in L^2(\RR^d), \nabla (hu)\in L^2(\RR^d), \forall \; h \in
\mathcal{C}^\infty_0(\RR^d)\},$$
 \be \label{W3}
 W_3(\RR^3):= \bigg\{ w \in H^1_{\mathrm{loc}}(\RR^3)  :
 \frac{w}{r}
 \in L^2(\RR^3), \nabla w \in
 L^2(\RR^3) \bigg\}
 \ee
and
 \be \label{W2}
 W_2(\RR^2) := \bigg \{ w \in H^1_{\mathrm{loc}}(\RR^2) :
 \frac{w}{ \sqrt{1+ r^2} \ln (2+ r^2)}  \in L^2(\RR^2), \nabla w
 \in L^2(\RR^2) \bigg\},
 \ee
where $r=|x|$. We also recall that $\Delta$ sets an isomorphism
from $W_d(\RR^d)$ to its dual $(W_d(\RR^d)^*$; see, for example,
\cite{nedelec}.

Observe that it is equivalent to the existence and uniqueness of
the solution in $W_d(\RR^d)$ to the problem (see, for instance,
\cite[Theorem 2.17]{book1})
\begin{equation}\label{eq:101-11}
\left\{
\begin{array}{ll}
\nabla \cdot\Gs \nabla u=0 & \mbox{in} \ \ B, \\
\GD u =0 & \mbox{in} \ \ \RR^d \setminus \overline{B}, \\
\ds \Gs \frac{\p u}{\p \nu} \Big|_- - \frac{\p u}{\p \nu} \Big|_+=g & \mbox{on} \ \ \p B,\\
u|_- - u|_+=f & \mbox{on} \ \ \partial B,\\
u(x) = O(|x|^{1-d})  & \mbox{as } |x| \rightarrow \infty.
\end{array}
\right.
\end{equation}
The injectivity of $\Acal$ comes directly from the uniqueness of a
solution to \eqnref{eq:101-1}. Since $u$ is
harmonic in $\RR^d \setminus \overline{B}$ and $u(x)
=O(|x|^{1-d})$ as $|x| \to \infty$, there exists $\phi \in
L^2_0(\p B)$ such that \be\label{uScalB} u(x) = \Scal_B [\phi]
(x), \quad x \in \RR^d \setminus \overline{B}. \ee If we set
$\psi=\sigma \pd{u}{\nu}|_-$, then
\begin{align}
u|_- = \Lambda_\sigma[\psi] +C,
\end{align}
where $C=\frac{1}{|\p B|} \int_{\p B} u|_-$. Note that
 \be\label{CScalB}
 C= \frac{1}{|\p B|} \int_{\p B} (u|_+ +f) = \frac{1}{|\p B|} \int_{\p B} \Scal_B [\phi].
 \ee
We now have from \eqnref{uScalB} and \eqnref{CScalB} that
\be\label{gpsi} g=\psi - \left(\frac{1}{2}I +
\mathcal{K}_B^*\right)[\phi] .\ee Furthermore, we have
\begin{align}
f=\Lambda_\sigma[\psi] + C - \Scal_B[\phi] = \Lambda_\sigma[\psi]-\tScal_B[\phi].
\end{align}
Thus $(\phi, \psi)$ satisfies \eqnref{Acalphi} and the proof is complete. \qed

We can now define the GPTs associated with $\Gs$ using the
operator $\Acal$.
\begin{definition} \label{def31}
Let $\Gs$ be a bounded measurable function in $\RR^d$, $d=2,3$,
such that $\Gs-1$ is compactly supported and \eqnref{glonetwo}
holds and let $B$ be a smooth domain satisfying \eqnref{suppB}.
For a multi-index $\alpha\in \NN^d$ with $|\alpha|\geq 1$, let
$(\phi_{\alpha},\psi_{\alpha}) \in H^{-1/2}(\p B) \times
H_0^{-1/2}(\p B)$ be the solution to
\begin{equation}\label{eq:104}
 \Acal \begin{bmatrix} \phi_\Ga \\ \psi_\Ga \end{bmatrix} = \begin{bmatrix} x^{\alpha}-
 \frac{1}{|\p B|} \int_{\p B} x^{\alpha} \;
 d\Ggs \\
 \nm
 \nu \cdot \nabla x^{\alpha} \end{bmatrix}
 \ \ \mbox{\rm on} \ \ \p B.
\end{equation}
For another multi-index $\beta \in \NN^d$, define the generalized
polarization tensors associated with the conductivity distribution
$\Gs(x)$ by \be\label{eq:GPTsdef}
M_{\alpha\beta}=M_{\alpha\beta}(\Gs) = \int_{\p B} x^{\beta}
\phi_{\alpha}(x) \, d\Ggs. \ee
\end{definition}

Definition \ref {def31} of the GPTs involves the domain $B$
satisfying \eqnref{suppB}. However, we will show later that GPTs for $\Gs$ (in fact, their harmonic combinations) are
independent of the choice of $B$ satisfying \eqnref{suppB}.

When $|\alpha|= |\beta|=1$, we denote $\BM:=(M_{\Ga\Gb})_{|\Ga|=|\Gb|=1}$ and call it the polarization tensor (matrix). Sometimes we write $\BM=(M_{ij})_{i,j=1}^d$.

For a given harmonic function $h$ in $\RR^d$, let $(\phi, \psi)$
be the solution to \eqnref{Acalphi2}. Since $$h(x)= h(0)+
\sum_{|\Ga| \ge 1} \frac{\p^\Ga h(0)}{\Ga !} x^\Ga,$$ we have
 \be
 \begin{bmatrix} \phi \\ \psi \end{bmatrix} = \sum_{|\Ga| \ge 1} \frac{\p^\Ga h(0)}{\Ga !} \begin{bmatrix} \phi_\Ga \\ \psi_\Ga \end{bmatrix}.
 \ee
By \eqnref{solrep} the solution $u$ to \eqnref{eq:101} can be written as
 $$
 u(x)=h(x) + \sum_{|\Ga| \ge 1} \frac{\p^\Ga h(0)}{\Ga !} \Scal_{B}[\phi_\Ga](x), \quad x \in \RR^d \setminus B.
 $$
Using the Taylor expansion
 $$
 \Gamma(x - y) = \sum_{|\Gb|=0}^{+\infty} \frac{(-1)^{|\Gb|}}{\Gb !} \p^\Gb \Gamma(x)  y^\Gb
 $$
which holds for all $x$ such that $|x| \to \infty$ while $y$ is
bounded \cite{book2}, we obtain the following theorem.
\begin{theorem}
For a given harmonic function $h$ in $\RR^d$,  let $u$ be the
solution to \eqnref{eq:101}. The following asymptotic formula
holds uniformly as $|x| \to \infty$:
 \be\label{uhasymp}
 u(x)-h(x) = \sum_{|\Ga|, |\Gb| \ge 1} \frac{(-1)^{|\Gb|}}{\Ga ! \Gb !} \p^\Ga h(0) M_{\Ga\Gb} \p^\Gb \Gamma(x).
 \ee
\end{theorem}

There is yet another way to represent the solution to
\eqnref{Acalphi2}. To explain it, let $\Lambda^e$ be the
NtD map
for the exterior problem:
$$ \Lambda^e[g]:= u|_{\p B}-\frac{1}{|\p B|} \int_{\p B} u,$$
where
$u$ is the solution to
\begin{equation}
\left\{
\begin{array}{ll}
\GD u =0 & \mbox{in }  \RR^d \setminus \overline{B}, \\
\ds \frac{\p u}{\p \nu} \Big|_+=g & \mbox{on } \p B,\\
u(x) =O(|x|^{1-d})  & \mbox{as } |x| \rightarrow \infty.
\end{array}
\right.
\end{equation}
Let $(\phi, \psi)$ be the solution to \eqnref{Acalphi2}. By \eqnref{gpsi}, we have
 \be\label{psiphih}
 \psi = \left(\frac{1}{2} I + \Kcal_B^*\right)[\phi] + \pd{h}{\nu}\Big|_{\p B}=\phi+\left(-\frac{1}{2} I + \Kcal_B^*\right)[\phi] + \pd{h}{\nu}\Big|_{\p B}.
 \ee
 On one hand, we obtain from the second identity in
\eqnref{psiphih} that $\int_{\p B} \phi =0$. On the other hand,
the first identity in \eqnref{psiphih} says that
 \be\label{psipdnu}
 \psi = \pd{}{\nu} \Scal_B [\phi] \Big|_{+} + \pd{h}{\nu}\Big|_{\p B} \quad\mbox{on } \p B,
 \ee
and hence
 \be\label{lamepsi}
 \Lambda^e[\psi]= \Scal_B[\phi] -\frac{1}{|\p B|} \int_{\p B}\Scal_B[\phi]  + \Lambda^e \left[\pd{h}{\nu}\Big|_{\p B}\right].
 \ee
Moreover,
 \be
 \psi = \phi + \pd{}{\nu} \Scal_B [\phi] \Big|_{-} + \pd{h}{\nu}\Big|_{\p B} \quad\mbox{on } \p B,
 \ee
and therefore,
 \be\label{lamonepsi}
 \Lambda_1[\psi]=\Lambda_1[\phi] + \Scal_B[\phi] + h|_{\p B} -\frac{1}{|\p B|} (\Scal_B[\phi]+ h).
 \ee
Combining \eqnref{lamepsi} and \eqnref{lamonepsi} with
 $$
 \Lambda_\sigma[\psi] = \Scal_B [\phi] + h|_{\p B}-\frac{1}{|\p B|} (\Scal_B[\phi]+ h)
 $$
in \eqnref{Acalphi2} yields
 \begin{align*}
 &(\Lambda_\sigma - \Lambda^e ) [\psi] = (\Lambda_1 - \Lambda^e)\left[ \pd{h}{\nu} \right],\\
 &(\Lambda_1 - \Lambda_\sigma)[\psi] = \Lambda_1[\phi].
 \end{align*}
Thus we readily get
 \begin{align}
 &\phi=\Lambda_1^{-1} (\Lambda_1 - \Lambda_\sigma) (\Lambda_\sigma - \Lambda^e )^{-1}
  (\Lambda_1 - \Lambda^e ) \left[ \pd{h}{\nu} \Big|_{\p B} \right],\label{eq:NtDeq}\\
 &\psi=(\Lambda_\sigma - \Lambda^e )^{-1} (\Lambda_1 - \Lambda^e ) \left[ \pd{h}{\nu} \Big|_{\p B}
 \right]. \label{eq:NtDeq2}
  \end{align}
Note that by the uniqueness  of a solution to problem
\eqnref{eq:101-11}, it is easy to see that $(\Lambda_\sigma -
\Lambda^e ):   H^{-1/2}_0(\partial B) \rightarrow
H^{1/2}_0(\partial B)$ is invertible.

Using \eqnref{eq:NtDeq} gives a slightly different but equivalent
definition of the GPTs.

\begin{lemma}
For all $\alpha, \beta \in \NN^d$,  $ M_{\alpha\beta}$, defined by
\eqnref{eq:GPTsdef}, can be rewritten in the following form:
\begin{equation} \label{getgpt} M_{\alpha\beta}(\Gs) = \int_{\p B} x^{\beta} \Lambda_1^{-1}
(\Lambda_1 - \Lambda_\sigma) (\Lambda_\sigma - \Lambda^e )^{-1}
  (\Lambda_1 - \Lambda^e ) \left[ \pd{x^\alpha}{\nu} \Big|_{\p B} \right] \,
d\Ggs. \end{equation}
\end{lemma}
Formula \eqnref{getgpt} shows how to get the GPTs from the NtD
maps.

%%%%%%%%%%%%%%%%%%%%%%%%%%%%%%%%%%%%%%%%%%%%%%%%%%%
\section{Properties of GPTs}\label{section4}
%%%%%%%%%%%%%%%%%%%%%%%%%%%%%%%%%%%%%%%%%%%%%%%%%%%

In this section, we prove important properties for the GPTs. We
emphasize that the harmonic sums of GPTs, not individual ones,
play a key role. Let $I$ and $J$ be finite index sets. Harmonic
sums of GPTs are $\sum_{\alpha \in I, \beta \in J} a_{\alpha}
b_{\beta} M_{\alpha\beta}$ where $\sum_{\alpha \in I} a_{\alpha}
x^{\alpha}$ and $\sum_{\beta \in J} b_{\beta} x^{\beta}$ are
harmonic polynomials.

The following lemma will be useful later.
\begin{lemma}\label{lemprop}  Let $I$ and $J$ be finite index sets.
Let $h_1(x):= \sum_{\alpha \in I} a_{\alpha} x^{\alpha}$ and
$h_2(x):= \sum_{\beta \in J} b_{\beta} x^{\beta}$ harmonic
polynomials and let $u_1$ be the solution to \eqnref{eq:101} with
$h_1(x)$ in the place of $h(x)$. Then, \be\label{sumab}
\sum_{\alpha \in I}\sum_{\beta \in J} a_{\alpha} b_{\beta}
M_{\alpha\beta} =\int_{\RR^d} (\sigma -1) \nabla u_1 \cdot \nabla
h_2 \, dx . \ee
\end{lemma}

\pf Let $\psi=\sum_{\alpha \in I}a_{\alpha} \psi_{\alpha}$ and $\phi=\sum_{\alpha \in I}a_{\alpha} \phi_{\alpha}$.
Then $u_1$ is given by
\begin{align*}
u_1(x):=\begin{cases}
h_1(x) + \Scal_B [\phi](x), & \ \ x \in \RR^d \setminus B, \\
\nm
 \Ncal_{B, \Gs} [\psi](x) + C, &  \ \ x \in B,
\end{cases}
\end{align*}

By \eqnref{eq:104}, \eqnref{eq:GPTsdef}, and the integration by parts, we see
\begin{align*}
\sum_{\alpha \in I}\sum_{\beta \in J} a_{\alpha} b_{\beta} M_{\alpha\beta}&= \int_{\p B} h_2(x) \phi(x) d\Ggs_x \\
&=\int_{\p B} h_2 \left( \pd{\mathcal{S}_B[\phi]}{\nu}\Big|_+  - \pd{\mathcal{S}_B[\phi]}{\nu}\Big|_- \right) d\Ggs_x \\
&=\int_{\p B} h_2 \left( \psi - \pd{h_1}{\nu} \right) d\Ggs_x -\int_{\p B} h_2 \pd{\mathcal{S}_B[\phi]}{\nu}\Big|_- d\Ggs_x \\
&=\int_{\p B} h_2 \left( \psi - \pd{h_1}{\nu} \right) d\Ggs_x -\int_{\p B} \mathcal{S}_B[\phi] \pd{h_2}{\nu} d\Ggs_x \\
&=\int_{\p B} h_2 \left( \psi - \pd{h_1}{\nu} \right) d\Ggs_x -\int_{\p B} \left( \GL_{\Gs}[\psi]- h_1 \right) \pd{h_2}{\nu} d\Ggs_x \\
&=\int_{\p B} \left( h_2 \psi - \GL_{\Gs}[\psi]\pd{h_2}{\nu} \right)  d\Ggs_x \\
&=\int_{\p B} \left( h_2 \sigma\pd{u}{\nu}\Big|_- - u\pd{h_2}{\nu} \right) d\Ggs_x \\
&=\int_{B} (\sigma -1) \nabla h_2 \cdot \nabla u_1 \, dx,
\end{align*}
which concludes the proof. \qed

Identity \eqnref{sumab} shows in particular that the definition of
(harmonic combinations of) the GPTs given in the previous section is independent of the
choice of $B$.

%%%%%%%%%%%%%%%%%%%%%%%%%%%%%%%%%%%%%%%%%%%%%%%%%%%%%%%%%%%%%%

%%%%%%%%%%%%%%%%%%%%%%%%%%%%%%%%%%%%%%%%%%%%%%%%%%%
\subsection{Symmetry}
%%%%%%%%%%%%%%%%%%%%%%%%%%%%%%%%%%%%%%%%%%%%%%%%%%%

We now prove symmetry of GPTs.

\begin{lemma}\label{th:sym} Let $I$ and $J$ be finite index sets.
For any harmonic coefficients $\{a_{\alpha}|\alpha \in I\}$ and
$\{b_{\beta}|\beta \in J\}$,  we have
\begin{equation}\label{eq:107}
\sum_{\alpha \in I}\sum_{\beta \in J} a_{\alpha} b_{\beta} M_{\alpha\beta} = \sum_{\alpha \in I}\sum_{\beta \in J} a_{\alpha} b_{\beta} M_{\beta\alpha}.
\end{equation}
In particular, the first-order GPT, $\BM$,  is symmetric.
\end{lemma}

\pf The symmetry property \eqnref{eq:107} can be easily deduced
from the proof of Lemma  \ref{lemprop}. However, we give here a
slightly different proof. For doing so, let
 $$
 h_1(x):= \sum_{\alpha \in I} a_{\alpha} x^{\alpha}, \quad
 h_2(x):= \sum_{\beta \in J} b_{\beta} x^{\beta}.
 $$
By \eqnref{eq:NtDeq}, we have
 \begin{align*}
\sum_{\alpha \in I}\sum_{\beta \in J} a_{\alpha} b_{\beta}
M_{\alpha\beta}
&= \int_{\p B} h_2(x) \phi(x) d\Ggs_x \\
& =  \int_{\p B} h_2 \, \Lambda_1^{-1} (\Lambda_1 - \Lambda_\sigma) (\Lambda_\sigma - \Lambda^e )^{-1} (\Lambda_1 - \Lambda^e ) \left[ \pd{h_1}{\nu} \Big|_{\p B} \right] \\
& =  \int_{\p B} \pd{h_2}{\nu} \, (\Lambda_1 - \Lambda_\sigma) (\Lambda_\sigma - \Lambda^e )^{-1} (\Lambda_1 - \Lambda^e ) \left[ \pd{h_1}{\nu} \Big|_{\p B} \right] .
\end{align*}
Since
$$ (\Lambda_\sigma -\Lambda^e)^{-1}=(\Lambda_1 -\Lambda^e)^{-1} + (\Lambda_\sigma -\Lambda^e)^{-1}(\Lambda_1-\Lambda_\sigma)(\Lambda_1 -\Lambda^e)^{-1},$$
we have
\begin{align*}
&\sum_{\alpha \in I}\sum_{\beta \in J} a_{\alpha} b_{\beta} M_{\alpha,\beta}\\
&= \int_{\p B} \pd{h_2}{\nu} \, (\Lambda_1 - \Lambda_\sigma ) \left[ \pd{h_1}{\nu} \right]+
\int_{\p B} \pd{h_2}{\nu}  \, (\Lambda_1 - \Lambda_\sigma) (\Lambda_\sigma - \Lambda^e )^{-1} (\Lambda_1 - \Lambda_\sigma ) \left[ \pd{h_1}{\nu} \right].
\end{align*}
Since the NtD maps,  $\Lambda_1, \Lambda_\sigma$, and $\Lambda^e$,
are self-adjoint, we get \eqnref{eq:107} which concludes the
proof. \qed

\subsection{Positivity and bounds}
%%%%%%%%%%%%%%%%%%%%%%%%%%%%%%%%%%%%%%%%%%%%%%%%%%%
Let $h(x)=\sum_{\alpha \in I} a_{\alpha} x^{\alpha}$ be a harmonic function in $\mathbb{R}^d$ and $u$ be the solution to \eqnref{eq:101}. As in the proof of Lemma \ref{lemprop}, we have
 \begin{align*}
\sum_{\alpha, \beta \in I} a_{\alpha} a_{\beta} M_{\alpha\beta}
 &=  \int_{\p B}  \left( h \, \sigma\pd{u}{\nu}\Big|_- - u\pd{h}{\nu} \right)ds\\
 &= \int_{\p B} \left(u\, \sigma \pd{u}{\nu}|_-  - 2(u-h) \pd{h}{\nu} -h \pd{h}{\nu}- (u-h)\pd{(u-h)}{\nu} \Big|_+ \right) ds\\
 &= \int_{B} \left(\sigma|\nabla u|^2  -2 \nabla(u-h)\cdot \nabla h - |\nabla h|^2\right) + \int_{\mathbb{R}^d\setminus\overline{B}} |\nabla(u-h)|^2\\
 &=\int_{\mathbb{R}^d}\left( \sigma |\nabla(u-h)|^2 +2(\sigma - 1) \nabla(u-h)\cdot \nabla h + (\sigma -1)
  |\nabla h|^2 \right)\\
 &=\int_{\mathbb{R}^d}\sigma\left| \nabla(u-h) +(1-\sigma^{-1})\nabla h\right|^2 + \int_B
 \frac{(\sigma -1)}{\sigma} |\nabla h|^2.
\end{align*}
We can also check the following variational principle:
\begin{equation}
\label{eq:109} \sum_{\alpha, \beta \in I} a_{\alpha} a_{\beta}
M_{\alpha\beta}= \min_{w\in W_d(\mathbb{R}^d)}\int_{\RR^d} \sigma
|\nabla w + (1-\Gs^{-1})\nabla h|^2+ \int_B \frac{(\sigma
-1)}{\sigma} |\nabla h|^2,
\end{equation}
where $W_d(\RR^d)$ is defined by (\ref{W3}) and (\ref{W2}).

Following the same lines of proof as in \cite{book2} for the
homogeneous case, we have the following bounds for GPTs.
\begin{theorem}
Let $I$ be a finite index set. Let $\{a_{\alpha} | \alpha \in I\}$
be the set of coefficients such that $h(x):=\sum_{\alpha \in I}
a_{\alpha} x^{\alpha}$ is a harmonic function. Then we have
\begin{equation}
\label{eq:111} \int_B\frac{(\Gs-1)}{\Gs} |\nabla h|^2 \leq
\sum_{\alpha,\beta \in I} a_{\alpha}a_{\beta} M_{\alpha\beta} \leq
\int_B (\Gs-1) |\nabla h|^2 .
\end{equation}
\end{theorem}
\pf The bound on the left-hand side is obvious since
$$\min_{w\in W_d(\mathbb{R}^d)}\int_{\RR^d} \sigma |\nabla w + (1-\Gs^{-1})\nabla h|^2\geq 0.$$
By taking $w=0$, we get
$$
\sum_{\alpha, \beta \in I} a_{\alpha} a_{\beta} M_{\alpha\beta}  \leq \int_B \frac{(\sigma -1)^2}{\sigma}
 |\nabla h|^2 + \int_B \frac{(\sigma -1)}{\sigma} |\nabla h|^2=\int_B (\Gs-1)|\nabla
 h|^2,
$$
which concludes the proof. \qed

The above theorem shows that if $\Gs$ is strictly lager than 1
then the GPTs are positive definite, and they are negative
definite if $0<\Gs<1$. Note that optimal bounds on the first-order
GPT have been derived in \cite{yvesvog,lipton}.

\subsection{GPTs and contracted GPTs}
%%%%%%%%%%%%%%%%%%%%%%%%%%%%%%%%%%%%%%%%%%%%%%%%%%%%%%%%%%%%%%

The contracted GPTs appeared in the asymptotic expansions as in
\eqnref{expan2} and \eqnref{expan3} while the GPTs appeared in
\eqnref{uhasymp}. By comparing those asymptotic formulas, we
obtain the following lemma which relates both quantities.
\begin{lemma} \label{lemgptvcontract}
\begin{itemize}
\item[\textrm{(i)}] If $r^n\cos n\theta = \sum_{|\alpha|=n}
a_{\alpha}^c x^\alpha$ and $r^n\sin n\theta = \sum_{|\alpha|=n}
a_{\alpha}^s x^\alpha$ in two dimensions, then
  \begin{align*}
M_{mn}^{cc}= \sum_{|\alpha|=m,|\beta|=n}   a_{\alpha}^c
a_{\beta}^c M_{\alpha\beta},\quad
M_{mn}^{cs}= \sum_{|\alpha|=m,|\beta|=n}   a_{\alpha}^c a_{\beta}^s M_{\alpha\beta},\\
M_{mn}^{sc}= \sum_{|\alpha|=m,|\beta|=n}   a_{\alpha}^s
a_{\beta}^c M_{\alpha\beta},\quad M_{mn}^{ss} =
\sum_{|\alpha|=m,|\beta|=n}   a_{\alpha}^s a_{\beta}^s
M_{\alpha\beta}.
 \end{align*}
\item[\textrm{(ii)}] If $r^n Y_n^m(\theta, \varphi)
 = \sum_{|\alpha|=n} a_{\alpha}^{mn} x^\alpha$ in three dimensions,
then
  \begin{align*}
 M_{mnkl}= \sum_{|\alpha|=n,|\beta|=l}  a_{\alpha}^{mn} a_{\beta}^{kl} M_{\alpha\beta}.
 \end{align*}
\end{itemize}
Conversely, any harmonic combination of the GPTs can be recovered
from the contracted GPTs.
\end{lemma}

%%%%%%%%%%%%%%%%%%%%%%%%%%%%%%%%%%%%%%%%%%%%%%%%%%%%%%%%%%%%%%
\subsection{Determination of NtD map}
%%%%%%%%%%%%%%%%%%%%%%%%%%%%%%%%%%%%%%%%%%%%%%%%%%%%%%%%%%%%%%

It is proved in \cite{AK02} (see also \cite[Theorem 4.9]{book2})
that the full set of harmonic combinations of GPTs associated with
a homogeneous inclusion determines the NtD map on the boundary of
any domain enclosing the inclusion, and hence the inclusion. In
the case of inhomogeneous conductivity inclusions, the same proof
can be easily adapted to obtain the following result.
\begin{theorem}\label{th:uni1}
Let $I$ and $J$ be finite index sets. Let $\Gs_i$, $i=1,2$, be two
 conductivity distributions with $\mbox{supp}~ (\sigma_i -1)
\subset \overline{B}$ and satisfying \eqnref{glonetwo}. If
 \be
 \sum_{\alpha \in I}\sum_{\beta \in J} a_{\alpha} b_{\beta} M_{\alpha\beta}(\Gs_1)= \sum_{\alpha \in I}\sum_{\beta \in J} a_{\alpha} b_{\beta} M_{\alpha\beta}(\Gs_2)
 \ee
for any harmonic coefficients $a_{\alpha}$ and $b_{\beta}$, then
 \be\label{NtDuniq}
 \GL_{\Gs_1}=\GL_{\Gs_2} \quad\mbox{on } \p B.
 \ee
\end{theorem}

Using uniqueness results of the Calder\'on problems (for example
\cite{SU87, AP06}) one can deduce from \eqnref{NtDuniq} that
$\Gs_1=\Gs_2$ under some regularity assumptions on the
conductivities imposed in those results. In two dimensions,
uniqueness holds for conductivities in $L^\infty$ \cite{AP06}.

%%%%%%%%%%%%%%%%%%%%%%%%%%%%%%%%%%%%%%%%%%%%%%%%%%%
\section{Sensitivity analysis for GPTs}
%%%%%%%%%%%%%%%%%%%%%%%%%%%%%%%%%%%%%%%%%%%%%%%%%%%

We now consider the sensitivity of the GPTs with respect to
changes in the conductivity distribution. Again, we suppose that
$\Gs-1$ is compactly supported in a domain $B$. The perturbation
of the conductivity $\Gs$ is given by $\Gs+\ep \gamma$, where
$\ep$ is a small positive parameter, $\gamma$ is compactly
supported in $B$ and refers to the direction of the changes.  The
aim of this section is to derive an asymptotic formula, as $\Ge
\to 0$, for the perturbation \be \triangle_{M}:=\sum_{\a \in I}
\sum_{\beta \in J}a_{\a}b_{\beta} \left(M_{\a\beta}(\Gs+\ep
\gamma)-M_{\a\beta}(\Gs)\right), \ee where $\{ a_\Ga \}$ and $\{
b_\Gb \}$ are harmonic coefficients and $I$ and $J$ are finite
index sets.

Let, as above, $h_1$ and $h_2$ be the harmonic functions given by
$$ h_1(x) = \sum_{\alpha} a_\alpha x^\alpha,\quad h_2(x) = \sum_{\beta} b_\beta x^\beta.$$
By \eqnref{eq:NtDeq} and a direct calculation we obtain
\begin{align*}
\triangle_{M}
&= \int_{\p B} h_2 \, \GL_1^{-1} \Big( (\GL_1- \GL_{\Gs+\ep\gamma})(\GL_{\Gs+\ep\gamma}-\GL^{e})^{-1}
 - (\GL_1- \GL_{\Gs})(\GL_{\Gs}-\GL^{e})^{-1} \Big)(\GL_1-\GL^e)[\na {h_1} \cdot \nu] \\
&= \int_{\p B} h_2 \, \GL_1^{-1}\Big( (\GL_{\Gs}- \GL_{\Gs+\ep\gamma})(\GL_{\Gs+\ep\gamma}-\GL^{e})^{-1} \\
& \quad +
(\GL_1- \GL_{\Gs})\left[(\GL_{\Gs+\ep\gamma}-\GL^{e})^{-1}
-(\GL_{\Gs}-\GL^{e})^{-1}\right]\Big)(\GL_1-\GL^e)[\na {h_1} \cdot \nu] \\
&= \int_{\p B} h_2 \, \GL_1^{-1} (\GL_{1}-\GL^e)(\GL_{\Gs}-\GL^e)^{-1} (\GL_{\Gs}-\GL_{\Gs+\ep\gamma})(\GL_{\Gs+\ep\gamma}-\GL^e)^{-1}
(\GL_1-\GL^e)[\na {h_1} \cdot \nu] \\
&= \int_{\p B} (\GL_{\Gs}-\GL_{\Gs+\ep\gamma})[g_2] \, g_1^{\ep},
\end{align*}
where
 \be\label{gonetwo}
 g_1^{\ep}=(\GL_{\Gs+\ep\gamma}-\GL^e)^{-1}(\GL_1-\GL^e)[\na h_1 \cdot \nu] \quad\mbox{and}\quad
 g_2=(\GL_{\Gs}-\GL^e)^{-1}(\GL_{1}-\GL^e)[\na h_2 \cdot \nu].
 \ee
Since $\GL_{\Gs}$ is self-adjoint, we have
\begin{align*}
\triangle_{M} & = \frac{1}{2} \int_{\p B} (\GL_{\Gs}-\GL_{\Gs+\ep\gamma})[g_2+g_1^{\ep}] \, (g_2+g_1^{\ep})
 - \frac{1}{2} \int_{\p B} (\GL_{\Gs}-\GL_{\Gs+\ep\gamma})[g_2] \, g_2\\
& \quad - \frac{1}{2}\int_{\p B} (\GL_{\Gs}-\GL_{\Gs+\ep\gamma})[g_1^{\ep}] \, g_1^{\ep}.
\end{align*}

We need the following two lemmas.
\begin{lemma}\label{le:NtD_eq}
If $u_1$ and $u_2$ are the solutions of $\nabla\cdot(\Gs_1\na
u_1)=0$ and $\nabla\cdot(\Gs_2\na u_2)=0$ with the Neumann
boundary conditions $\Gs_1\frac{\p u_1}{\p \nu}=g$ and
$\Gs_2\frac{\p u_2}{\p \nu}=g$ on $\p B$, respectively, then the
following identity holds \be\label{eq:NtD_eq} \int_{\p B}
(\GL_{\Gs_2} - \GL_{\Gs_1})[g] \, g \,d\Ggs = \frac{1}{2}\int_{B}
(\Gs_1-\Gs_2) \Big(|\na (u_1-u_2)|^2 + |\na u_1|^2 +|\na
u_2|^2\Big) dx. \ee
\end{lemma}
\pf The following identity is well-known (see, for instance,
\cite{book2}):
$$
\int_{B} \Gs_1 |\na (u_1-u_2)|^2 dx + \int_{B}(\Gs_1-\Gs_2)|\na
u_1|^2 dx = \int_{\p B} (\GL_{\Gs_2} - \GL_{\Gs_1})[g] \, g
\,d\Ggs .
$$
We also have
$$
\int_{B} \Gs_2 |\na (u_1-u_2)|^2 dx - \int_{B}(\Gs_1-\Gs_2)|\na
u_2|^2 dx = \int_{\p B} (\GL_{\Gs_1} - \GL_{\Gs_2})[g] \, g
\,d\Ggs.
$$
Subtracting those two equalities we obtain (\ref{eq:NtD_eq}). \qed

\begin{lemma}
\label{le:est_op}
There is a constant $C$ such that
\be\label{DtNstability}
\|\GL_{\Gs_1}-\GL_{\Gs_2}\|  \leq  C \|\Gs_1-\Gs_2\|_{L^{\infty}(B)} .
\ee
\end{lemma}
\pf Let $u_1$ and $u_2$ be the solutions to $\nabla\cdot(\Gs_1\na u_1)=0$ and $\nabla\cdot(\Gs_2\na u_2)=0$ with boundary conditions $\Gs_1\frac{\p u_1}{\p \nu}=g$ and $\Gs_2\frac{\p u_2}{\p \nu}=g$ on $\p B$, respectively, and let
$h$ be the harmonic function with $\pd{h}{\nu} = g$ on $\p B$. Then we have
$$
\int_{B} \Gs_1 |\na u_1|^2 = \int_{B} \na h \cdot \na u_1 \leq \frac{1}{2\epsilon} \int_{B} |\na h|^2 dx+ \frac{\epsilon}{2}\int_{B} |\na u_1|^2 dx,
$$
for any $\epsilon>0$. Choosing $\epsilon=\inf_{B} \Gs_1:=\underline{\Gs_1}$ we get
$$
\int_{B} |\na u_1|^2 dx\leq \frac{1}{\underline{\Gs_1}^2}\int_{B} |\na h|^2dx.
$$
Similarly, we get
$$
\int_{B} |\na u_2|^2 dx\leq \frac{1}{\underline{\Gs_2}^2}\int_{B}
|\na h|^2dx,
$$
where  $\underline{\Gs_2}: = \inf_{B} \Gs_2$. It then follows from \eqnref{eq:NtD_eq} that
\begin{align*}
\left|\int_{\p B} (\GL_{\Gs_2} - \GL_{\Gs_1})[g] \, g
\,d\Ggs\right| & \leq  \frac{3}{2}\|\Gs_1-\Gs_2\|_{L^\infty(B)}
\left( \int_{B}|\na u_1|^2 dx+ \int_{B}|\na u_2|^2 dx \right) \\
& \leq \frac{3}{2} (\frac{1}{\underline{\Gs_1}^2}+\frac{1}{\underline{\Gs_2}^2}) \big(\int_{B} |\na h|^2dx \big)\, \|\Gs_1-\Gs_2\|_{L^\infty(B)} \\
& \leq C \| g \|_{H^{-1/2}(\p B)}^2 \, \|\Gs_1-\Gs_2\|_{L^\infty(B)}.
\end{align*}
Thus, we obtain \eqnref{DtNstability}. \qed

With the notation \eqnref{gonetwo} in hand, let
 \be\label{gone}
 g_1:= (\GL_{\Gs}-\GL^e)^{-1}(\GL_{1}-\GL^e)[\na h_1 \cdot \nu],
 \ee
and let $u_i$, for $i=1,2,$ be the solution to \be\label{uonetwo}
\left\{
\begin{array}{ll}
\na \cdot (\Gs \na u_i)=0 & \mbox{in }  B, \\
\ds \Gs \frac{\p u_i}{\p \nu}=g_i & \mbox{on } \p B.
\end{array}
\right.
\ee
Let $u_1^\ep$, $v_1^\ep$ and $v_2^\ep$ be the solutions to
\be
\left\{
\begin{array}{ll}
\na \cdot (\Gs \na u_1^\ep)=0 & \mbox{in }  B, \\
\nm \ds \Gs \frac{\p u_1^\ep}{\p \nu}=g_1^\ep & \mbox{on } \p B,
\end{array}
\right.
\ee
\be
\left\{
\begin{array}{ll}
\na \cdot \big( (\Gs  + \ep \gamma)  \na v_1^\ep \big)=0 & \mbox{in }  B, \\
\nm \ds \Gs \frac{\p v_1^\ep}{\p \nu}=g_1^\ep & \mbox{on } \p B,
\end{array}
\right.
\ee
and
\be
\left\{
\begin{array}{ll}
\na \cdot \big( (\Gs  + \ep \gamma) \na v_2^\ep \big)=0 & \mbox{in }  B, \\
\ds \Gs \frac{\p v_2^\ep}{\p \nu}=g_2 & \mbox{on } \p B.
\end{array}
\right. \ee Then by Lemma \ref{le:NtD_eq} we have
\begin{align*}
\triangle_{M} & =  \frac{1}{4}\int_{B} \ep\gamma\Big(|\na (u_2 + u_1^\ep - v_2^{\ep} - v_1^{\ep})|^2 + |\na (u_2 + u_1^\ep)|^2 + |\na (v_2^{\ep} + v_1^{\ep})|^2\Big) dx \\
& \quad - \frac{1}{4}\int_{B} \ep\gamma\Big(|\na (u_2 - v_2^{\ep})|^2+ |\na u_2|^2 +|\na v_2^{\ep}|^2\Big) dx  \\
& \quad - \frac{1}{4}\int_{B} \ep\gamma\Big(|\na (u_1^\ep - v_1^{\ep})|^2+ |\na u_1^\ep|^2 + |\na v_1^{\ep}|^2\Big) dx .
\end{align*}
Lemma \ref{le:est_op} yields
 \be
 \| u_1^\ep - u_1 \|_{L^2(B)}^2 = O(\ep)
 \ee
and
 \be
 \| v_j^\ep - u_j \|_{L^2(B)}^2 = O(\ep), \quad j=1,2.
 \ee
Thus we get
 $$
 \triangle_{M} = \frac{\ep}{2}\int_{B} \gamma \Big(|\na (u_1+u_2)|^2- |\na u_1|^2-|\na u_2|^2\Big) dx + O(\ep^2) ,
 $$
to arrive at the following theorem.
\begin{theorem}\label{th:pertub}
Let $I$ and $J$ be finite index sets. Let $u_1$ and $u_2$ be the
solutions to \eqnref{uonetwo}. Then we have \be\label{eq:pertub}
\sum_{\a \in I} \sum_{\beta \in J} a_{\a}b_{\beta}
M_{\a\beta}(\Gs+\ep\gamma)= \sum_{\a \in I} \sum_{\beta \in J}
a_{\a}b_{\beta} M_{\a\beta}(\Gs)+\ep \int_{B} \gamma \na u_1 \cdot
\na u_2 dx + O(\ep^2). \ee
\end{theorem}

%%%%%%%%%%%%%%%%%%%%%%%%%%%%%%%%%%%%%%%%%%%%%%%%%%%
\section{Reconstruction of an inhomogeneous conductivity distribution}
%%%%%%%%%%%%%%%%%%%%%%%%%%%%%%%%%%%%%%%%%%%%%%%%%%%
Over the last decades, a considerable amount of work has been
dedicated to the inverse conductivity problem. We refer, for
instance, to \cite{borcea, cheney} and the references therein.

Here, our approach is completely different. We stably recover some
important features of inhomogeneous conductivities using their
GPTs. It should be emphasized that the GPTs can be obtained from
boundary measurements by solving a least-squares problem
\cite{han}. The purpose of this section is to illustrate
numerically the viability of this finding.

For doing so, we use a least-square approach (see, for instance,
\cite{hanke}). Let $\Gs^*$ be the exact (target) conductivity (in
two dimensions) and let $y_{mn}:=M_{mn}(\Gs^*)$ (omitting for the
sake of simplicity $c$ and $s$ for the superscripts in contracted
GPTs). The general approach is to minimize over bounded
conductivities $\Gs$ the discrepancy functional \be\label{eq:lsq}
S(\Gs)=\frac{1}{2}\sum_{m+n \le N}
\omega_{mn}\|y_{mn}-M_{mn}(\Gs)\|^2 \ee for some finite number $N$
and some well-chosen weights $\omega_{mn}$. The weights
$\omega_{mn}$ are used to enhance resolved features of the
conductivity as done in \cite{siims, borcea2}. We solve the above
minimization problem using the gradient descent (Landweber)
method.

%%%%%%%%%%%%%%%%%%%%%%%%%%%%%%%%%%%%%%%%%%%%%%%%%%%
\subsection{Fr\'echet derivative and an optimization procedure}
%%%%%%%%%%%%%%%%%%%%%%%%%%%%%%%%%%%%%%%%%%%%%%%%%%%

Let, again for the sake of simplicity,  $M_{mn}(\sigma)=
M_{mn}^{cc}(\sigma)$ be the contracted GPTs for a given
conductivity $\Gs$. The Fr\'echet derivative in the direction of
$\Gg$, $M'_{mn}(\Gs)[\Gg]$, is defined to be
 $$
 M'_{mn}(\Gs)[\Gg] := \lim_{\ep \to 0} \frac{M_{mn}(\sigma+ \ep \Gg)-M_{mn}(\sigma)}{\ep}.
 $$
From (\ref{eq:pertub}) we obtain  that \be\label{eq:Fre_op}
M'_{mn}(\Gs)[\gamma]= \int_{B} \gamma \na u_n \cdot \na u_m dx,
\ee where $u_n$ and $u_m$ are the solutions of
\begin{equation} \label{defh}
\left\{
\begin{array}{ll}
\na \cdot (\Gs \na u)=0 & \mbox{in }  B, \\
\nm \ds \Gs \frac{\p u}{\p
\nu}=(\GL_{\Gs}-\GL^e)^{-1}(\GL_{1}-\GL^e)[\na h \cdot \nu] &
\mbox{on } \p B,
\end{array}
\right.
\end{equation}
with $h= r^n \cos n\theta$ and $h=r^m \cos m\theta$, respectively.
Note that if $M_{mn}(\sigma)$ is one of the other contracted GPTs,
then $h$ should be changed accordingly.

One can easily see that the adjoint $M'_{mn}(\Gs)^*$ of
$M'_{mn}(\Gs)$ is given by
 \be
 M'_{mn}(\Gs)^*[c]=c \na u_m \cdot \na u_n, \quad c \in \RR.
 \ee
The gradient descent procedure to solve the least-square problem
\eqnref{eq:lsq} reads \be \label{eq:Land1} \Gs_{k+1}= \Gs_k +
\sum_{m,n}\omega_{mn} M'_{mn}(\Gs_k)^*[y_{mn}-M_{mn}(\Gs_k)] .\ee
In the numerical implementation, the GPTs for the exact
conductivity distribution can be computed by using the following
formula: \be \label{eq:solM} M_{mn}(\Gs)= \int_{\p B}
(h_1-\tilde{u}_1) \Gs \frac{\p u_2}{\p \nu} d\Ggs= \int_{\p B}
h_1\Gs \frac{\p u_2}{\p \nu} d\Ggs- \int_{\p B} \frac{\p h_1}{\p
\nu} u_2 d\Ggs ,\ee where $\tilde{u}_1$ and $u_2$ are the
solutions to \be\label{eq:solm} \left\{
\begin{array}{ll}
\nabla \cdot(\Gs(x)\nabla \tilde{u}_1)=0 & \mbox{in} \ \ B, \\
\nm \ds \Gs\frac{\partial \tilde{u}_1}{\partial \nu} = \frac{\p
h_1}{\p \nu} & \mbox{on} \ds \ \ \partial B \ \ (\int_{\p B}
\tilde{u}_1=0),
\end{array}
\right.
\ee
and
\begin{equation}\label{eq:soln}
\left\{
\begin{array}{ll}
\nabla \cdot \Gs \nabla u_2 =0 & \mbox{in} \ \ \RR^d, \\
\nm u_2(x)-h_2(x)=O(|x|^{1-d}) & |x| \rightarrow \infty,
\end{array}
\right.
\end{equation}
respectively. Here, $h_1= r^n \cos n \theta$ and $h_2= r^m \cos m
\theta$ in two dimensions.

% and \be M_{m,n}'(\Gs)^*%\na u_m \cdot \na u_n . \ee We see from (\ref{eq:solM}) that we
%only need to solve (\ref{eq:soln}) to get the GPTs and its
%conjugate gradient.

%However, in the numerical implementation, we can get more accurate
%$\tilde{u}_m$ by solving (\ref{eq:solm}) than $u_n$ by solving
%(\ref{eq:soln}).

%Note that, in real applications, the GPTs can be obtained from the
%NtD map.

%we can not solve (\ref{eq:solm}) and (\ref{eq:soln}). But we can
%measure $\tilde{u}_m$ and calculate the inverse of
%$\GL_{\Gs}-\GL_{1}^e$ (see \eqnref{eq:NtDeq}) in order to get the
%GPTs and its conjugate gradient.

On the other hand, in order to compute $M'_{mn}(\Gs)^*$ we need to
invert the operator $\GL_{\Gs}-\GL^e$. This can be done
iteratively. In fact,  the least-square solution to
$$
(\GL_{\Gs}-\GL^e)[g]=f,
$$
is given by \be \label{eq:Land2} g_{k+1}=g_k +
\omega(\GL_{\Gs}-\GL^e)(f-(\GL_{\Gs}-\GL^e)[g_k]),\ee where
$\omega$ is a positive step-size.

In order to  stably and accurately reconstruct the conductivity
distribution, we use a recursive approach proposed in
\cite{AKLZ12} (see also \cite{siims, AGKLY11, borcea2, bao}). We
first minimize the discrepancy between the first contracted GPTs
for $1\leq m, n \leq l$. Then we use the result as an initial
guess for the minimization between the GPTs for $1\leq m, n \leq l
+1$. This corresponds to choosing appropriately the weights
$\omega_{mn}$ in \eqnref{eq:lsq}. Moreover, we refine the mesh
used to compute the reconstructed conductivity distribution every
time we increase the number of used contracted GPTs in the
discrepancy functional.

\subsection{Resolution analysis in the linearized case}
Let $d=2$ and let $B$ be a disk centered at the origin. Consider
the linearized case by assuming that the conductivity $\sigma$ is
given by $\sigma = k + \ep \gamma$, where $k \neq 1$ is a positive
constant and $\ep$ is a small parameter. In that case, using
Theorem \ref{th:pertub} together with Lemma \ref{lemgptvcontract},
one can easily see that
$$
M_{mn}(k + \ep \gamma) = M_{mn}(k) + \ep \int_B \gamma \nabla u_m
\cdot \nabla u_n \, dx + O(\ep^2),
$$
with $u_m(x)=r^m e^{i m\theta}, u_n(x) = r^n e^{i n \theta},$ and
$x=(r,\theta)$. Hence, it follows that
$$\bigg( \int_B \gamma(r,\theta) r^{m+n-2} e^{i (m+n)\theta}\,
d\theta \, dr \bigg)_{1\leq m,n \leq N}$$ can be obtained from the
contracted GPTs, $M_{mn}$, for $1\leq m,n\leq N$. Therefore, the
higher is $N$, the better is the angular resolution in
reconstructing $\gamma$. On the other hand, it is clear that
variations of $\gamma$ that are orthogonal (in the $L^2$ sense) to
the set of polynomials $(r^{m+n-2})_{1\leq m,n\leq N}$ cannot be
reconstructed from the contracted GPTs $M_{mn}$, $1\leq m,n\leq
N$. Moreover, the reconstruction of $\gamma$ near the origin
($r=0$) is more sensitive to noise than near the boundary of $B$.
This is in accordance with \cite{procams, gunther}.

%%%%%%%%%%%%%%%%%%%%%%%%%%%%%%%%%%%%%%%%%%%%%%%%%%%
\subsection{Numerical illustration}
%%%%%%%%%%%%%%%%%%%%%%%%%%%%%%%%%%%%%%%%%%%%%%%%%%%
In this section, for simplicity we only consider the
reconstruction from contracted GPTs of a conductivity distribution
which is radially symmetric. Many recent works have been devoted
to the reconstruction of radially symmetric conductivities. See,
for instance, \cite{radia1, radia2, radia3}.

Here we consider the following conductivity distribution:
 \be
 \Gs=(0.3r^2+0.5r^3+6(r^2-0.5)^2+3.0)/3.0,
 \ee
and apply our
original approach for recovering $\Gs$ from the contracted GPTs
$M_{mn}$, for $m,n\leq N$. Since the conductivity distribution
$\Gs$ is radially symmetric we have
\begin{eqnarray*}
M_{mn}^{cs}&=&M_{mn}^{sc}=0 \ \ \mbox{for} \ \ \mbox{all} \  \ m,n, \\
M_{mn}^{cc}&=&M_{mn}^{ss}=0 \ \ \mbox{if} \ \ m\ne n,
\end{eqnarray*}
and $M_m:=M_{mm}^{cc}=M_{mm}^{ss}$. We use $M_1$ to estimate the
constant conductivity which has the same first-order GPT as
follows:
 \be
 \sigma_0:= \frac{2|B|+ M_1}{2|B|- M_1}.
 \ee
Then we use $\sigma_0$ as an initial guess and apply the recursive
approach described below.

\begin{figure}[h]
\begin{center}
  \includegraphics[width=4.5in,height=3.3in]{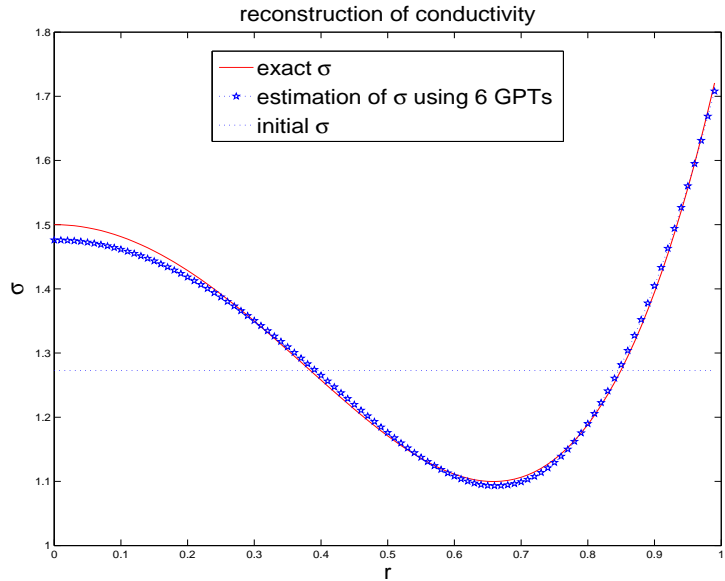}
  \end{center}
  \caption{Reconstructed conductivity distribution. \label{fig-result}}
\end{figure}

\begin{figure}[h]
\begin{center}
  \includegraphics[width=4.5in,height=3.3in]{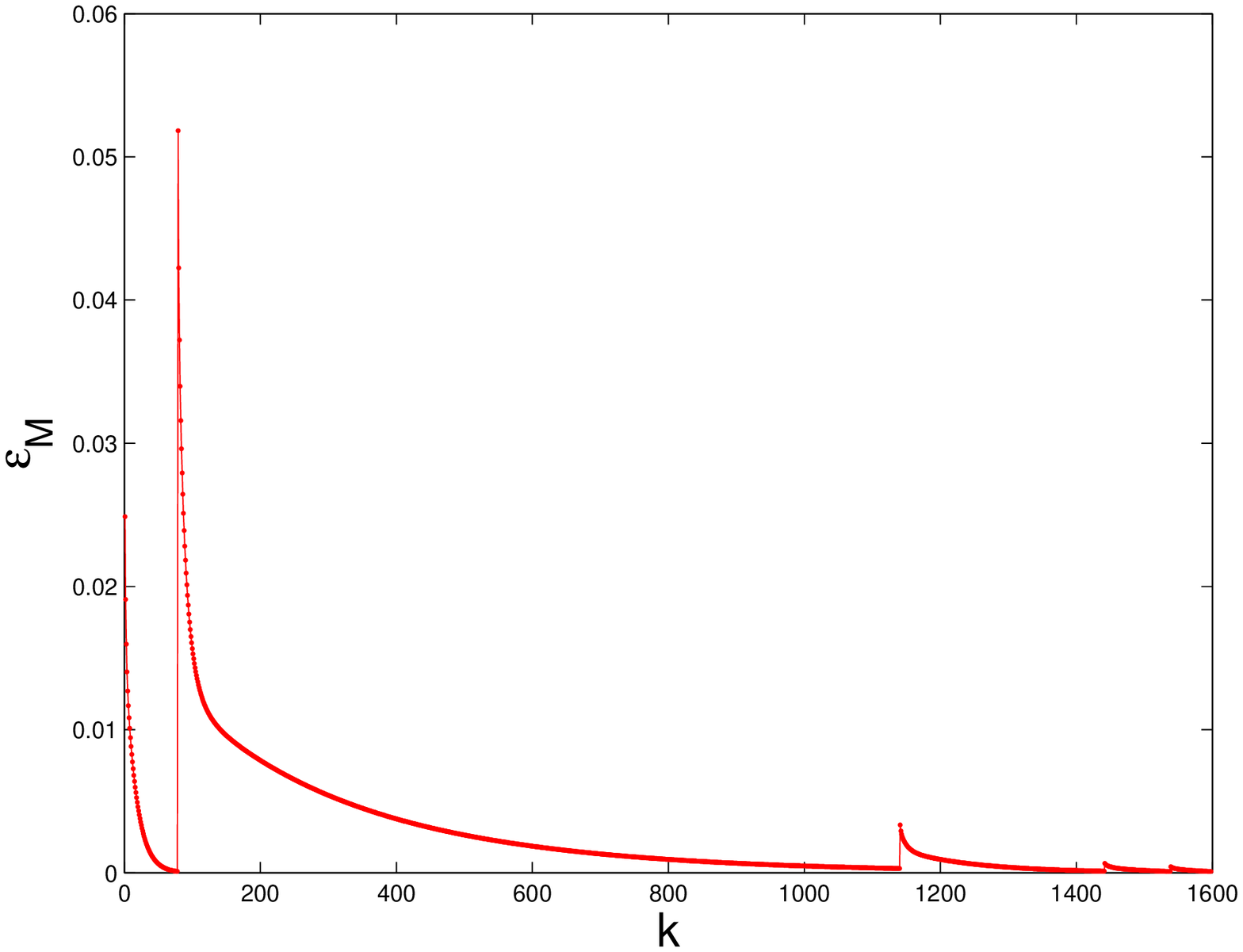}
  \end{center}
  \caption{The convergence history of $\varepsilon_{M}$, where ${\rm k}$ is the number of iterations. \label{fig-result1}}
\end{figure}

\begin{figure}[h]
\begin{center}
  \includegraphics[width=4.5in,height=3.3in]{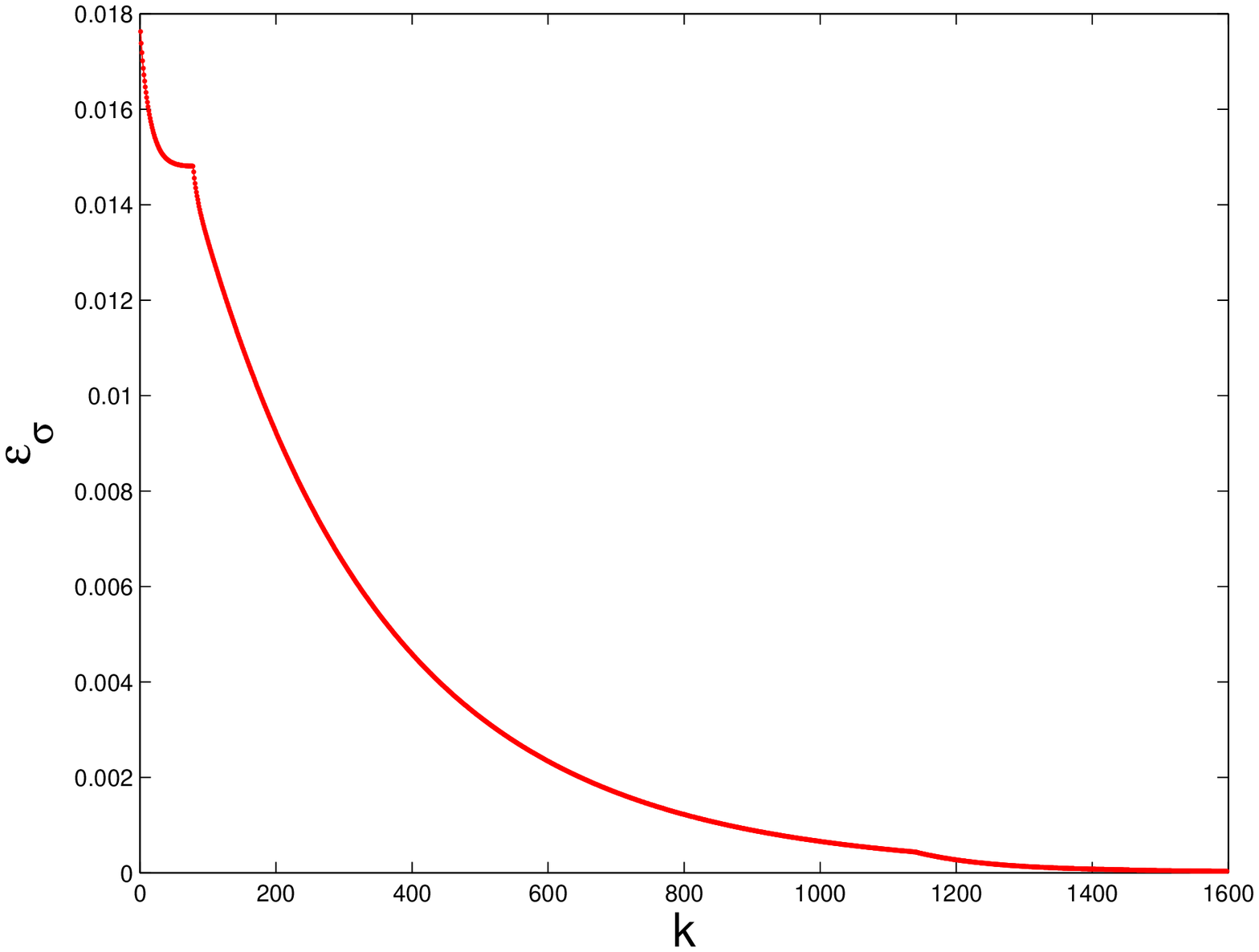}
  \end{center}
  \caption{The convergence history of $\varepsilon_{\Gs}$, where ${\rm k}$ is the number of iterations. \label{fig-result2}}
\end{figure}

Let $k_*$ be the last iteration step, and let $\varepsilon_{M}$ and $\varepsilon_{\Gs}$ be discrepancies of GPTs and the conductivities, {\it i.e.},
 \be
 \varepsilon_{M}:= \sum_{n \le N}(y_{n}-M_{n}(\Gs_{k_*}))^2, \quad y_n:= M_{n}(\Gs),
 \ee
($N$ represents the number of GPTs used) and
 \be
 \varepsilon_{\Gs}:=\frac{\int_{ B}
(\Gs_{k_*}-\Gs)^2}{\int_{ B} \Gs^2}.
 \ee
Figure \ref{fig-result} shows the reconstructed conductivity distribution
using contracted GPTs with $N=6$. In this reconstruction, the errors $\varepsilon_{M}$ and
$\varepsilon_{\Gs}$ are given by
$$\varepsilon_{M}=9.83318e-005, \quad
\varepsilon_{\Gs}=3.5043e-005,$$ after $1598$ iterations. It
should be noted that the conductivity is better reconstructed near
the boundary of the inclusion than inside the inclusion itself.
Figure \ref{fig-result1} shows how fast $\varepsilon_{M}$
decreases as the iteration proceeds. The sudden jump in the figure
happens when we switch the number of GPTs from $N$ to $N+1$.
Figure \ref{fig-result2} is for the convergence history of
$\varepsilon_{\Gs}$.

\section{Conclusion}
In this paper we have introduced for the first time the notion of
GPTs for inhomogeneous conductivity inclusions. The GPTs carry out
overall properties of the conductivity distribution. They can be
determined from the NtD map. We have established positivity and
symmetry properties for the GPTs. We have also analyzed their
sensitivity with respect to small changes in the conductivity. We
have proposed a recursive algorithm for reconstructing the
conductivity from the GPTs and presented a numerical example
to show that radially symmetric conductivities can be accurately
reconstructed from the GPTs. A numerical study of the use of the
GPTs for solving the inverse conductivity problem will be the
subject of a forthcoming work. A stability and resolution analysis
will be performed. It would also be very interesting to extend the
ideas of this paper to the inverse wave medium problems.

\end{document}